\documentclass[11pt,a4paper]{article}
\usepackage[utf8]{inputenc}
\usepackage[T2A]{fontenc}
\usepackage[english]{babel}
\usepackage{amsmath,amssymb,amsthm}
\usepackage{mathtools}
\usepackage{hyperref}
\usepackage{subcaption}
\usepackage{tikz}
\usepackage{pgfplots}
\pgfplotsset{compat=1.18}
\usepackage{geometry}
\usepackage{enumitem}
\usepackage{setspace}
\usepackage{comment}
\usepackage{graphicx}
\hypersetup{
    colorlinks=true,
    linkcolor=blue,
    citecolor=blue,
    urlcolor=blue
}
\usepackage{titleps}
\newpagestyle{main}{%
\sethead[\thepage][][] 
        {Spectral structure of infinite size squared distances matrices}{}{\thepage} 
}
\title{\textbf{Spectral structure of infinite size squared distances matrices}}

\author{\textsc{Alexander Plakhotnikov \footnote{
\textsc{Student. Saint-Petersburg State University, Department of Mathematics and Mechanics.} \\
\textit{Email adress:} \texttt{st132770@student.spbu.ru}}}}
\date{}
 
\begin{document}

\maketitle

{\small 
\textsc{Abstract.} Let a finite set of points $\{\xi_1,...,\xi_k\}$ be chosen in a metric space $(X,d)$, and let the squared distance matrix $\mathfrak{D}=(\mathfrak{D}(\xi_i,\xi_j)^2)_{i,j=1}^{k}$ be constructed from them. We propose a geometric approach to studying the spectral properties of squared distance matrices of infinite size, constructed from a countable set of points $\{\xi_k\}_{k\in \mathbb{Z}}$ on Riemannian manifold $(M,g)$. We move from the discrete problem to a continuous one using walk matrices. We describe the structure of the spectrum and study the properties of spectral flows. 
}

\bigskip

\begin{center}
    1. \textsc{Introduction}
\end{center}

Given a finite set of points $\xi_1,\dots,\xi_n$ in a metric space $(X,d)$, the squared distance matrix is defined by
\[
\mathfrak{D} = \big(d(\xi_i,\xi_j)^2\big)_{i,j=1}^n 
\]

The problem of studying the spectral properties of squared distance matrices reduces to studying a broader class of matrices called hollow symmetric nonnegative matrices~\cite{HSM}:

$$HSN(\mathbb{R}) = \bigg\{ A\in \operatorname{Mat}_{n\times n}(\mathbb{R}) \; \; | \;\;a_{ii}=0 \; , \; a_{ij}\geq0  \; \; \forall i\neq j\bigg\}$$

The possibility of zero off-diagonal elements is often excluded:

$$HSN^+(\mathbb{R}) = \bigg\{ A\in \operatorname{Mat}_{n\times n}(\mathbb{R}) \; \; | \;\;a_{ii}=0 \; , \; a_{ij}>0 \; \; \forall i\neq j \bigg\}$$

For the case of a finite set of points $\{\xi_1,...\xi_k\}$, the properties of the eigenvalues of $HSN(\mathbb{R})$ matrices are well known. A classical result, known as the Perron-Frobenius theorem, states that in a metric space $(X,\mathfrak{D})$ equipped with the Euclidean metric there exists a unique eigenvalue $\lambda^*\in \sigma(\mathfrak{D})$ such that 

$$\begin{cases}
    \lambda^*>0 \\
    \lambda^* = \max_{\lambda \in \sigma(\mathfrak{D})} \lambda
\end{cases}$$

To generalize the problem of studying the spectrum of a matrix of squared distances constructed from a countable set of points $\{\xi_k\}_{k\in \mathbb{Z}}$, it is necessary to introduce appropriate terminology. Following ~\cite{Tri} and ~\cite{Infd}, we introduce the definition of an infinite matrix as follows:

\medskip

\textbf{Definition 1.} \textit{A mapping $\mathfrak{M}: \mathbb{Z}\times\mathbb{Z} \longmapsto \mathbb{R}$ will be called an abstract infinite matrix (or simply an abstract matrix) with coefficients from $\mathbb{R}$.} For the elements of the matrix, we introduce the notation $\mathfrak{M}(n,m) \;\;n,m\in \mathbb{Z}$.

\medskip

\textbf{Definition 2.} \textit{The sequence $\alpha_k : \mathbb{Z} \longmapsto \mathbb{R}$, where $\alpha_k(n) = \mathfrak{M}(n,n-k), \; \; n\in \mathbb{Z}$ will be called the $k$-th diagonal of the matrix $\mathfrak{M}$. For $k=0$, we obtain the main diagonal}.

\medskip

\textit{The sequence $\beta_k : \mathbb{Z} \longmapsto \mathbb{R}$, $\beta_k(n) = \mathfrak{M}(k,n), \;\; n,k\in \mathbb{Z}$ is the $k$-th row of the abstract matrix $\mathfrak{M}$}.
\medskip

\textit{The sequence $\gamma_k : \mathbb{Z} \longmapsto \mathbb{R}$, $\gamma_k(n) = \mathfrak{M}(n,k), \;\; n,k\in \mathbb{Z}$ is the $k$-th column of the abstract matrix} $\mathfrak{M}$.

\medskip

\textbf{Definition 3.} \textit{Let $\mathfrak{A}_{n\times n}$ and $\mathfrak{B}_{n\times n}$, where $n = \infty$, be abstract matrices. We say that the product of abstract matrices is defined}

$$\mathfrak{A\cdot B = C} = \big(\mathfrak{C}(m,n)\big) = \bigg(\sum_{l=1}^\infty \mathfrak{A}(m,l)\cdot \mathfrak{B}(l,n)\bigg) \quad m,n \in \mathbb{Z}$$

\textit{if $c_{mn}$ is a convergent series} $\forall m,n$.
\medskip

\medskip

\textbf{Remark 1.} According to our definition, the product of abstract matrices is defined if the corresponding rows and columns of the factors form convergent series. Further, we will show that the rows and columns of the abstract matrix of squared distances form divergent series. Therefore, the product of the abstract matrix of squares with other abstract matrices is not defined. This circumstance makes spectral (and other) decomposition of abstract matrices of squared distances difficult.

\medskip

\textbf{Definition 4.} \textit{The principal minor of order $p\in \mathbb{N}$ of an abstract matrix $\mathfrak{M}$ will be called a finite square submatrix $D_p$ of size $p$, with elements} 

$$ d_{mk} = \mathfrak{M}(m,k) ; \;\;\;|m|,|k| \leq p $$

\medskip

Note that in the case of a finite-size matrix, the spectrum of squared distances depends on the metric.  

\medskip

\textbf{Definition 5.} \textit{By the spectrum $\sigma(\mathfrak{M})$ of an abstract matrix we will mean}:
\begin{enumerate}
    \item \textit{The union of the spectra of the principal minors of all orders of the abstract matrix. This set is denoted as}
    
    $$\sigma_{d}(\mathfrak{M}) = \bigcup_{p\in \mathbb{N}} \sigma(D_p)$$
    
    \textit{and will hereafter be called the discrete spectrum of the abstract matrix. The multiplicity of a number $\lambda \in \sigma_d(\mathfrak{M})$ is the number of its occurrences in the elements of the union.}
    
    \item \textit{If convergent subsequences $\{\lambda_{n_k}\}$ can be selected from the elements of the discrete spectrum $\sigma_d(\mathfrak{M})$, then we say that the limits of these subsequences $\{\mu_k\}$ form the continuous spectrum $\sigma_c(\mathfrak{M})$}
\end{enumerate}

Note that our definition of spectrа can be considered as local in this paper. In the case of an infinite-size matrix, the spectrum  depends on the metric. Further, we will show that the overall structure of the spectrum does not depend on the metric. From the last definition, it is clear that the spectrum is a subset of $\mathbb{R}$. Note that in the general case, the spectrum of an abstract matrix  $\sigma(\mathfrak{M})=\sigma_d(\mathfrak{M})\bigcup \sigma_c(\mathfrak{M})$ can be everywhere dense in some subset $S\subset \mathbb{R}$. Our definition of the spectra of an abstract matrix allows us to determine the transition to an abstract Jordan normal form, which can also be considered as local in this paper.\\

\textbf{Remark 2.} We define the abstract diagonal matrix as

$$\mathfrak{R}=\big(\mathfrak{R}(m,n) \big), \;  \;
\mathfrak{R}(m,n) =
\begin{cases}
    \mathfrak{R}(m), \quad m=n \\
    0, \quad\quad\quad m\neq n
\end{cases} \;  \quad\quad \forall m \;\;\;\mathfrak{R}(m)\in \mathbb{R} $$

\textbf{Definition 6.}  We say that the trace is defined for an abstract matrix if the sum of the eigenvalues from $\sigma_d(\mathfrak{M})$ is finite: 

$$\operatorname{Tr} \mathfrak{M} := \sum_{\lambda \in \sigma_d(\mathfrak{M})} \lambda < + \infty$$

This definition of the trace (taking into account our definition of the spectrum) allows, as we will show below, to see the structure of the spectrum without defining an abstract Jordan normal form. It is clear that such a series of eigenvalues will converge for abstract squared distances matrices, again due to our definition of the spectrum.

\begin{center}
2. \textsc{Main results}
\end{center}

Let a countable set of points $\{\xi_k\}_{k\in \mathbb{Z}}$ be arbitrarily chosen on Riemannian manifold $(M,g)$, and let the abstract squared distance matrix $\mathfrak{\mathfrak{D}}$ be constructed from them (the distances between points are geodesic):
\[
\mathfrak{\mathfrak{D}} (n,m)= \operatorname{dist}(\xi_n,\xi_m)^2 = d_{nm} \quad n,m \in \mathbb{Z}
\]

Now let the chosen points depend on time: $$  \quad\xi_k = \xi_k(t) \quad k \in \mathbb{Z}$$ 
The path described by the points (which is finite) will be naturally called a "walk", and we impose natural restrictions:

$$\begin{cases}
     \dim M >0 \\
     \sup_{x,y \in M}\operatorname{dist(x,y)}<+\infty
\end{cases}$$

Note, first, that the metric $g$ can be arbitrary. Second, we further work on the time interval $\mathcal{T}=[t_0;t_1]$. 
The squared distance matrix also becomes dependent on time $\mathfrak{D} = \mathfrak{D}(t)$. The concept of a \emph{walk matrix} arises:

$$ \mathcal{WM}(\mathbb{R}) = \bigg\{ \mathfrak{D}(t)  \; \bigg| \;\forall t \in \mathcal{T} \; d_{ii}(t)=0 \; , \; d_{ij}(t)>0  \; \; \; \forall i\neq j \bigg\}$$

\medskip

Clearly, the set $\mathcal{WM}(\mathbb{R})$ forms a convex cone without the origin. Note that homotopy of paths $[\xi_k]$ implies homotopy of walk matrices. On this set, there is an equivalence relation given by continuous matrix homotopy ($I=[0,1]$)
\[
\forall U(t),V(t) \in \mathcal{WM}(\mathbb{R}) \quad \exists \Phi :  I\times \mathcal{T}  \longmapsto \mathcal{WM}(\mathbb{R})
\]
\[
\Phi(0,t)= U(t) \quad  \Phi(1,t)= V(t) 
\]

\begin{center}
\textsc{2.1. Recovering the matrix from the spectrum, structure of the spectrum}    
\end{center}

All statements carry over to $HSN^{+}(\mathbb{R})$ as a special case.

\medskip
 
\textbf{Lemma 1.} \textit{There is no sequence} $\{x_n\}_{n\in\mathbb{Z}}\subset \mathbb{R}$ \textit{converging to an element of} $\mathbb{R}\setminus \{0\}$ \textit{such that its terms are elements of the spectrum of the matrix} $\mathfrak{D}(t)\in \mathcal{WM}(\mathbb{R}) \; \; \forall t \in \mathcal{T} $.

\textbf{}

\textit{Proof.} If such a sequence exists, then due to the known properties of the squared matrix, it must satisfy the following condition: since $\forall t \in \mathcal{T} \; \operatorname{Tr}\mathfrak{D}(t)=0$, then $\sum_nx_n=0$. In addition, if such a sequence exists, then its elements belong to a discrete spectrum (by our definition), and the limit of this sequence belongs to a continuous spectrum (by our definition).  But since it must converge in $\mathbb{R}\setminus\{0\}$, then $x_n \longmapsto \mu \neq 0$. Contradiction with the condition $ \sum_nx_n=0$. $\square$

\medskip

\textbf{Lemma 2.}  $\forall t \in [t_0;t_1] \quad \forall \mathfrak{D}(t)\in \mathcal{WM}(\mathbb{R})$ \textit{the discrete spectrum can only accumulate at zero.}
    
\medskip

\textbf{Proof.} Proof by contradiction. Suppose the discrete spectrum accumulates at an arbitrary nonzero number $\mu \neq 0$. Then we get a contradiction with Lemma 1. $\square$

\medskip
From Lemma 2 we directly obtain

\medskip

\textbf{ Lemma 3. } $\forall t \in [t_0;t_1] \quad \forall \mathfrak{D}(t)\in \mathcal{WM}(\mathbb{R})$ \textit{The discrete spectrum cannot be everywhere dense in any connected (more than one-point) subset of} $\mathbb{R}$. 

\medskip

\textbf{Theorem 1.} $\forall t \in [t_0;t_1] \quad \forall \mathfrak{D}(t)\in \mathcal{WM}(\mathbb{R})$ \textit{only three cases of spectrum structure are possible}:
\medskip
\begin{enumerate}
    \item \textit{The spectrum is purely continuous. Then} $\sigma(\mathfrak{D}(t))=\{0\}$.
    \item \textit{The spectrum is purely discrete. Then it consists of a finite set of numbers.}
    \item \textit{The spectrum is the union of point and continuous spectra. Then the discrete spectrum accumulates at the point} $\{0\}$ -- \textit{the only point of the continuous spectrum.}
\end{enumerate} 
\medskip

\textit{Proof.} 
Note, first of all, that in the general case of an arbitrary self-adjoint operator in an infinite-dimensional separable Hilbert space, a continuous spectrum can consist of more than a single-point connected subset $\mathbb{R}$. Perhaps even from the entire real axis.
\medskip
\begin{enumerate}
    \item 
    \begin{enumerate}
        \item[1.1] The continuous spectrum cannot be a more than one-point connected subset of the negative or positive semi-axis, because such a set is uncountable, while the number of points $\lambda \in \sigma_c(\mathfrak{M})$ is at most countable by the definition.

        \item[1.2] Two options remain: either it is $\{0\}$, or a countable or finite set of numbers. Note that the point continuous spectrum, unlike the general case of operators, where the continuous spectrum consists in part of connected intervals (which are excluded as stated), arises as the limit of numbers from the discrete spectrum. Then Lemma 1 excludes the possibility of at least one nonzero number. Only zero remains, which is guaranteed by Lemma 2.
    \end{enumerate}

    \item Indeed, let us show that the case of at least a countable set is impossible. Suppose the spectrum is purely discrete and consists of a countable set of points. Then there are at most two options:
    \begin{itemize}
        \item[2.1] It accumulates at zero. But then we get a one-point continuous spectrum $\{0\}$. Contradiction with the condition of pure discreteness.
   
        \item[2.2] It accumulates at a nonzero number. This option is excluded by Lemma 1.
    \end{itemize}
    \item Follows from the previous items and Lemmas 1 and 2. $\square$
\end{enumerate}

\medskip

From what has been proven we obtain

\medskip

\textbf{Corollary.} \textit{Every sequence} $\{x_n\}_{n\in \mathbb{Z}}\subset \mathbb{R}$ \textit{converging to zero defines the spectrum of some squared distance matrix.}

\medskip
\textbf{Theorem 2.} \textit{For all walk matrices} $\forall t \in \mathcal{T} \quad \forall \mathfrak{D}(t)\in \mathcal{WM}(\mathbb{R})$ \textit{the series}
$$||\mathfrak{D}||_{\infty} = \max_i \sum_j |d_{ij}|$$
 \textit{diverges}.

\medskip

\textbf{Proof.} Proof by contradiction. Suppose $\exists j_0$ for which $\max_{j_0} \sum_i |d_{ij_0}| = C<+\infty$ is attained. Then, by the necessary condition for series convergence, we have $d_{ij_0}\longmapsto 0$ as $i\mapsto \infty$. But if the maximum is finite, then the other sums are also finite

$$\forall j \quad \sum_i |d_{ij}| \leq \max_{j_0} \sum_i |d_{ij_0}| = C< +\infty$$

Consequently, by the necessary condition for series convergence, 

$$\forall j \quad d_{ij} \longmapsto 0  \;\;  i\mapsto \infty$$

We obtain:
\begin{enumerate}
    
\item for $j=1$, the distance between the first point and the $i$-th point is arbitrarily small: $d_{i1} \longmapsto 0$ as $i\longmapsto \infty$. This means the sequence of selected points ${\{\xi_k\}}_{k\in \mathbb{Z}}$ accumulates at $\xi_1$.
\item for $j=2$, similarly $d_{i2} \longmapsto 0$ as $i\mapsto 0$. This means the sequence of points $\{\xi_k\}_{k\in \mathbb{Z}}$ accumulates at $\xi_2$.\\

    ...
    and so on.
\end{enumerate}
That is, $\forall j \in \mathbb{Z}$ the sequence of points ${\xi_k}$ accumulates at $\xi_j$. Which is impossible. $\square$
 
\medskip

The natural question arises about the connection between the spectrum and the norm of the $HSN^+(\mathbb{R})$. Let us consider the matrices  $\mathfrak{D} \in HSN^+(\mathbb{R})$ that define the operators $\mathfrak{D}:\ell_2 \longmapsto \ell_2$. Then we have

\medskip
\textbf{Theorem 3.} The spectrum $\sigma(\mathfrak{D})$ is nontrivial, i.e., contains nonzero elements $\iff$ $$||\mathfrak{D}||_{\infty} = \max_{i \in \mathbb{Z}} \sum_{j \in \mathbb{Z}} |d_{ij}| <+\infty$$
\textit{Proof.} 
\begin{itemize}
    \item[$\Longleftarrow$] Since the series with elements
\[||\mathfrak{D}||_{\infty} = \max_{i \in \mathbb{Z}} \sum_{j \in \mathbb{Z}} |d_{ij}| < +\infty\]
converges \(\implies ||\mathfrak{D}||_{\infty} < C \in \mathbb{R}\). Let's construct a sequence from the Weyl criterion. Let's show that for the matrix \(\mathfrak{D}\) one can construct a sequence \(\{x_n\}_{n \in \mathbb{Z}}\) such that
\[||x_n|| = 1, \quad ||\mathfrak{D} \cdot x_n - \lambda \cdot E \cdot x_n|| \to 0\]
Indeed, consider an arbitrary sequence for which $||x_n||=1 \;\;\forall n$:
\[||\mathfrak{D} \cdot x_n - \lambda \cdot E \cdot x_n||^2 = \left\langle \mathfrak{D} \cdot x_n - \lambda \cdot E \cdot x_n, \, \mathfrak{D} \cdot x_n - \lambda \cdot E \cdot x_n \right\rangle =\]
\[= ||Dx_n|| + \lambda^2||x_nE|| - 2\lambda \left\langle Dx_n, \, x_nE \right\rangle \leq\]
\[\leq ||\mathfrak{D}|| \cdot ||x_n|| + \lambda^2 ||E|| \cdot ||x_n|| + 2\lambda \sqrt{||\mathfrak{D}||} \cdot ||x_n||\]
\[= ||\mathfrak{D}|| + 2\lambda \sqrt{||\mathfrak{D}||} + \lambda^2 = \left( \sqrt{||\mathfrak{D}||} + \lambda \right)^2\]
Since, obviously, $||\mathfrak{D}||_{\infty} \neq 0$, we obtain
\[\lambda \longrightarrow -\sqrt{||\mathfrak{D}||} \neq 0 \Longrightarrow \exists \lambda \neq 0, \lambda \in \sigma(\mathfrak{D})\]
Hence follows the nontriviality of the spectrum \(\sigma(\mathfrak{D})\). Moreover, we have shown, without using the property $\operatorname{Tr}\mathfrak{D}=0$, that there exists at least one negative eigenvalue.

\item[$\Longrightarrow$] By condition, the spectrum is nontrivial. That is, \(\exists \lambda \in \sigma(\mathfrak{D}) : \lambda \neq 0\). According to the Weyl criterion for this \(\lambda\), there exists a sequence \(\{x_n\}_{n \in \mathbb{Z}}\), such that

\[||x_n|| = 1, \quad ||\mathfrak{D} \cdot x_n - \lambda \cdot E \cdot x_n|| \longrightarrow 0\]
 We obtain
\[0 \leftarrow ||\mathfrak{D} \cdot x_n - \lambda \cdot E \cdot x_n||^2 = ||Dx_n|| + \lambda^2||x_nE|| - 2\lambda \left\langle Dx_n, x_nE \right\rangle\]
And two cases are possible:
1) \(\lambda > 0\) Then

\[||Dx_n|| + \lambda^2||x_nE||  -  2\lambda\left\langle Dx_n, x_nE \right\rangle > -2\lambda \left\langle Dx_n, x_nE \right\rangle\]

But this tends to zero if and only if \(\lambda \longmapsto 0\). Contradiction with the assumption of nontrivial spectrum.

2) \(\lambda < 0\)  Then

$$ ||Dx_n|| + \lambda^2||x_nE||  -  2\lambda\left\langle Dx_n, x_nE \right\rangle > ||Dx_n|| + \lambda^2$$

the latter sum in $\mathbb{R}$ can never tend to zero. That is, we have obtained

\[0 \leftarrow ||\mathfrak{D} \cdot x_n - \lambda \cdot E \cdot x_n||^2 \geq ||Dx_n|| + \lambda^2 \not\longrightarrow 0\]

which is impossible and is an explicit contradiction with the Weyl criterion for the chosen \(\lambda\) $\square$
\end{itemize}

\medskip

\textbf{Remark.} It is easy to see that the series of eigenvalues (accounting for multiplicities) $\lambda\in \sigma (\mathfrak{D})$ will always converge. Indeed,  
$$0=\operatorname{Tr}\mathfrak{D}  = \sum_{\lambda \in \sigma_d(\mathfrak{D})}\lambda = \sum_{p}\operatorname{Tr}D_p$$

From where it follows immediately\\

\textbf{Corollary.} $\forall \varepsilon > 0$ interval $\mathbb{R}\setminus(-\varepsilon;\varepsilon)$ contains only a finite set of points of the $\sigma(\mathfrak{D})$.

\newpage

\begin{center}
\textsc{2.2. Property of the spectral flow of $\mathcal{WM}$--matrices}    
\end{center}

In the theory of operators with a discrete spectrum, for a continuous path of operators $\{B_t\}_{t\in [0;1]}$, the concept of spectral flow $\operatorname{sf}(\{B_t\}_{t\in[0;1]})$ is used, defined as the number of eigenvalues (taking into account multiplicities) passing through the point $0$ in the positive direction. The spectral structure of infinite size squared distances matrices allows the use of spectral flow as one of the methods of spectrum analysis.

Obviously, in the general case $\sigma(\mathfrak{D}(t))$ is not a homotopy invariant. Therefore, the spectral flow is also not a homotopy invariant.  Moreover, since the spectrum is completely determined by the location of points on a Riemannian manifold (with a fixed metric), the spectral flow (with an  appropriate choice of paths) is not limited above.\\
Let $\mathfrak{Diff}$ be the group of diffeomorphisms of the manifold $(M,g)$. 

\medskip

\textbf{Property.} In the general case, $\sigma(\mathfrak{D}(t))$ is not an invariant of the action \\ $\mathfrak{Diff} \curvearrowright (M,g)$.

\medskip

As proof, we give a simple example.\\
\textit{Example.} Let $(M,g) = (\mathbb{S}^n_R \; , \;g)$ be the $n$-dimensional sphere of radius $R$ with an arbitrary metric. Then for each $t^*\in \mathcal{T}$ consider the mapping

$$f : \mathbb{S}^n_R \longmapsto \mathbb{S}^n_{q(t^*) R}$$

\medskip

Stretching by a factor of $q(t^*)>1$. Clearly, increasing the pairwise distances between paths leads to an increase in all elements of the corresponding walk matrix $\mathcal{WM}(\mathbb{R})$ and, as a consequence, to a change in the eigenvalues $\lambda_i(t)$.

\medskip

It is also obvious that for a fixed manifold $M$ and metric $g$ the spectral flow is completely determined by the displacement of the points $\xi_i(t)$.
\medskip

\medskip

\textbf{Open problem.} In special cases, it is possible to establish a finite upper bound for the spectral flow. It is clear that on any manifold $M$ satisfying conditions $(*)$, by regulating the walk of the chosen points $\{\xi_i(t)\}$, one can also regulate the spectral flow. Therefore, in the general case, for a fixed time interval $[t_0;t_1]$, the upper estimate $0\leq |\operatorname{sf}| \leq C$ can also be changed. A formal proof of the independence of the finite upper estimate from the walk on specific manifolds is difficult.

\begin{center}
    \textsc{2.3. The growth rate of positive and negative eigenvalues }
\end{center}

An interesting observation is that the growth rate of the number of negative and positive eigenvalues depends on the metric. Below, as an example, are graphs for a three-dimensional basic space with the Minkowski metric 
$$\mathfrak{D}(\xi_1,\xi_2)= \bigg ( \sum_{i=1}^3 |\xi_1^i -\xi^i_2|^p\bigg )^{1/p}$$

  for different $p$. In our method of generating points in space, the distribution is formed from the sum of three components: a uniform distribution, a Student's t-distribution, and a centered and scaled gamma distribution. The question of the theoretical justification of the asymptotic increment of negative and positive eigenvalues for specific metrics remains open.

\medskip

\begin{figure}[h]
    \centering
    \includegraphics[width=0.45\textwidth]{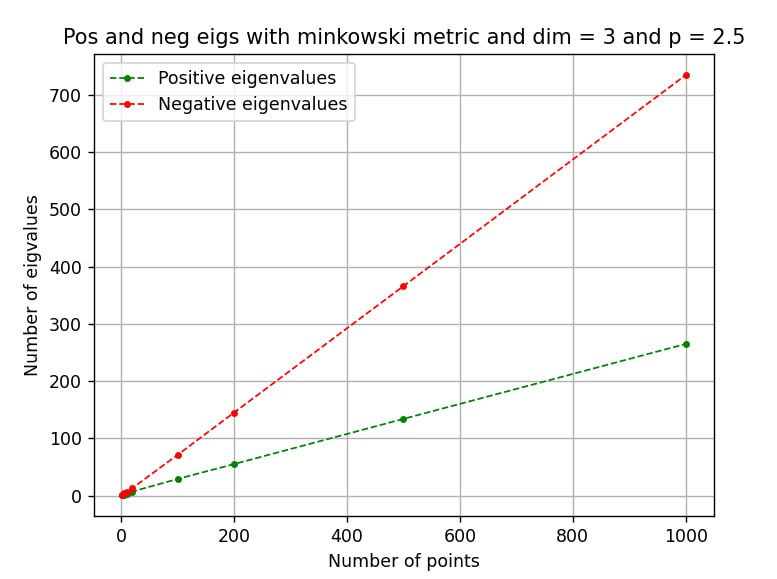}
    \hfill
    \includegraphics[width=0.45\textwidth]{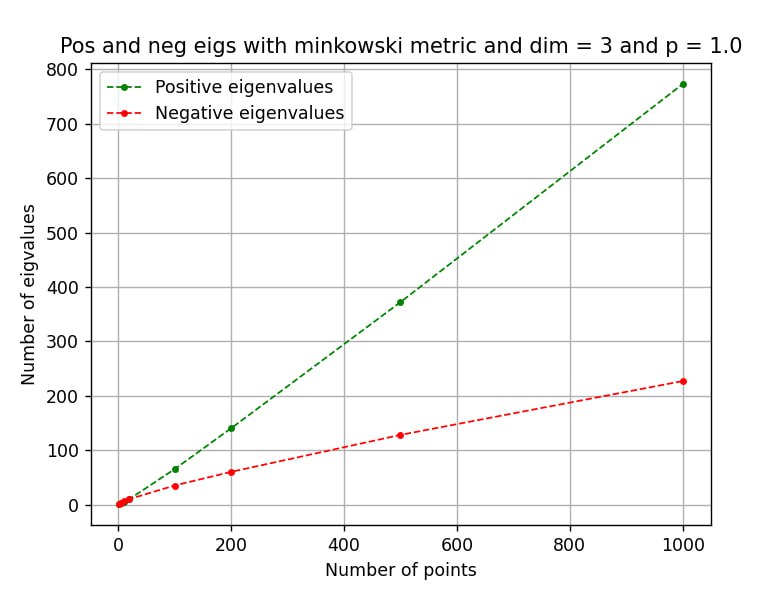}
\end{figure}

\bigskip

\end{document}